\documentclass{ifacconf}

\usepackage{graphicx}
\usepackage{tikz}
\usepackage{natbib}
\usepackage{amsmath}
\usepackage{amssymb}
\usepackage{mathtools}
\usepackage{url}

\DeclareMathOperator{\Real}{Re}
\DeclareMathOperator{\Imag}{Im}
\DeclareMathOperator{\sign}{sign}
\DeclareMathOperator{\diff}{d\!}
\DeclareMathOperator{\id}{I}

\DeclarePairedDelimiter{\abs}{\lvert}{\rvert}
\newenvironment{proof}{\begin{pf}}{\hspace*{\fill}$\square$\end{pf}}
\newenvironment{proof*}[1]{\begin{pf*}{#1}}{\hspace*{\fill}$\square$\end{pf*}}

\begin{document}
\begin{frontmatter}

\title{Insights into the multiplicity-induced-dominancy for scalar delay-differential equations with two delays}

\author[LSS,Inria]{Sébastien Fueyo}
\author[LSS,Inria]{Guilherme Mazanti}
\author[LSS,Inria,IPSA]{Islam Boussaada}
\author[LSS]{Yacine Chitour}
\author[LSS,Inria]{Silviu-Iulian Niculescu}

\address[LSS]{Universit\'e Paris-Saclay, CNRS, CentraleSup\'elec, Laboratoire des signaux et syst\`emes, 91190, Gif-sur-Yvette, France. \\ (e-mails: \{first name.last name\}@centralesupelec.fr)}
\address[Inria]{Inria, Saclay--Île-de-France Research Center, DISCO Team, France.}
\address[IPSA]{Institut Polytechnique des Sciences Avancées (IPSA), 63 boulevard de Brandebourg, 94200 Ivry-sur-Seine, France.}

\begin{abstract}
It has been observed in recent works that, for several classes of linear time-invariant time-delay systems of retarded or neutral type with a single delay, if a root of its characteristic equation attains its maximal multiplicity, then this root is the rightmost spectral value, and hence it determines the exponential behavior of the system, a property usually referred to as multiplicity-induced-dominancy (MID). In this paper, we investigate the MID property for one of the simplest cases of systems with two delays, a scalar delay-differential equation of first order with two delayed terms of order zero. We discuss the standard approach based on the argument principle for establishing the MID property for single-delay systems and some of its limitations in the case of our simple system with two delays, before proposing a technique based on crossing imaginary roots that allows to conclude that the MID property holds in our setting.
\end{abstract}


\end{frontmatter}

\section{Introduction}

Time delays appear naturally in several systems stemming from applications in engineering, physics, economics, chemistry, or biology, for instance, as they are useful to model in a simplified manner the effects of finite-speed propagation of mass, energy, or information. For that reason, the stability analysis, control, and stabilization of time-delay systems have been extensively studied in the scientific literature and were the subject of several monographs, such as \cite{1963differential, diekmann, gu2003stability, hale1993introduction, michiels2014stability, Sipahi2011Stability}.

A question of major interest on the stability of time-delay systems is whether, given some choice on the parameters of the system, they can be chosen in such a way as to render the system exponentially stable, sometimes with some additional constraints on performance criteria for the system. This question appears naturally, for instance, in the feedback stabilization of time-delay control systems, in which the aim is to choose the input of the system as a function of its state (possibly with delays) in order to render the resulting system exponentially stable.

Consider a linear time-delay system of retarded type and with finitely many pointwise delays under the general form
\begin{equation}
\label{eq:general-LTI-delay}
x^\prime(t) = \sum_{j=0}^N A_j x(t - \tau_j),
\end{equation}
where $N$ is a nonnegative integer, $x(t) \in \mathbb R^d$ is the (instantaneous) state of the system, $d$ is a positive integer, $\tau_0, \dotsc, \tau_N$ are nonnegative delays, and $A_0, \dotsc, A_N$ are $d \times d$ matrices with real coefficients. The \emph{spectrum} of \eqref{eq:general-LTI-delay} is the set of roots of its \emph{characteristic function} $\Delta: \mathbb C \to \mathbb C$ defined by
\begin{equation}
\label{eq:spectrum}
\Delta(s) = \det\left(s \id - \sum_{j=0}^N A_j e^{-s \tau_j}\right),
\end{equation}
and \eqref{eq:general-LTI-delay} is exponentially stable if and only if all roots of $\Delta$ have negative real part (as proved, for instance, in Chapter~7 of \cite{hale1993introduction}). Clearly, this is equivalent to requiring that the \emph{rightmost root} of $\Delta$ in $\mathbb C$ has negative real part.

The function $\Delta$ from \eqref{eq:spectrum} is a \emph{quasipolynomial}, i.e., an entire function $Q$ that can be written under the form $Q(s) = P_0(s) + \sum_{j = 1}^M e^{s r_j} P_j(s)$, where $P_0, \dotsc, P_M$ are polynomials and $r_1, \dotsc, r_M$ are pairwise distinct nonzero real numbers. A classical result provided in \cite[Part Three, Problem~206.2]{Polya1998Problems} and known as \emph{P\'olya--Szeg\H{o} bound} implies that the multiplicity of any root of $Q$ is at most equal to the \emph{degree} $D$ of $Q$, defined as the integer $D = M + \sum_{j=0}^M d_j$, where $d_0, \dotsc, d_M$ are the degrees of the polynomials $P_0, \dotsc, P_M$, respectively.

From the spectral point of view, the stabilization problem for \eqref{eq:general-LTI-delay} amounts to asking, if one has the choice of some coefficients of (some of) the matrices $A_1, \dotsc, A_N$, whether they can be chosen in order for all the roots of $\Delta$ to have negative real part. Except in the trivial case where all delays are equal to zero, $\Delta$ always admits infinitely many roots. However, since $A_1, \dotsc, A_N$ have finitely many coefficients, one cannot assign arbitrarily the roots of $\Delta$ in the complex plane. A natural idea is then to look for methods for assigning only finitely many roots of $\Delta$ but guaranteeing that the rightmost root of $\Delta$ is among those assigned. If such a method is available, then the stabilization problem can be solved by using that method to assign roots with negative real parts.

A possible such method, at least for some classes of time-delay systems, is to assign a real root of $\Delta$ with maximal multiplicity. Indeed, some works have shown that, for some classes of time-delay systems, a real root of maximal multiplicity is necessarily the rightmost root, a property known as \emph{multiplicity-induced-dominancy} (MID). This link between maximal multiplicity and dominance has been suggested in \cite{Pinney1958Ordinary} after the study of some simple, low-order cases, but without any attempt to address the general case and, up to the authors' knowledge, very few works have considered this question in more details until recently in works such as \cite{Boussaada2016Multiplicity, Boussaada_ESAIM, BOUSSAADA_further_remark, ramirez16, Mazanti2020Spectral, mazanti2021multiplicity, MBNC, Benarab2020MID}. These works consider only time-delay equations with a single delay and the MID property was shown to hold, for instance, for retarded equations of order $1$ in \cite{Boussaada2016Multiplicity}, proving dominance by introducing a factorization of $\Delta$ in terms of an integral expression when it admits a root of maximal multiplicity $2$; for retarded equations of order $2$ with a delayed term of order zero in \cite{BOUSSAADA_further_remark}, using also the same factorization technique; or for retarded equations of order $2$ with a delayed term of order $1$ in \cite{Boussaada_ESAIM}, using Cauchy's argument principle to prove dominance of the multiple root. Most of these results are actually particular cases of a more general result on the MID property from \cite{mazanti2021multiplicity} for retarded equations of order $n$ with delayed term of order $n-1$, which relies on links between quasipolynomials with a root of maximal multiplicity and Kummer's confluent hypergeometric function. The MID property was also extended to neutral systems of orders $1$ and $2$ in \cite{Ma,Benarab2020MID, MBNC}, as well as to the case of complex conjugate roots of maximal multiplicity in \cite{Mazanti2020Spectral}.

This paper addresses the MID property for one of the first nontrivial time-delay systems with more than one delay, namely the scalar delay-differential equation with two delays
\begin{equation}
\label{eq:main-syst}
y^\prime(t) + a_0 y(t) + a_1 y(t - \tau_1) + a_2 y(t - \tau_2) = 0,
\end{equation}
where $a_0, a_1, a_2 \in \mathbb R$ and $\tau_1, \tau_2 \in (0, +\infty)$ with $\tau_1 \neq \tau_2$. With no loss of generality, we assume in this paper that $0 < \tau_1 < \tau_2$. In addition to the fact that several real-world systems naturally have more than one delay, an important motivation for studying the stability and stabilization of systems with several delays is the fact that additional delays may, in some situations, improve stability properties, as shown, for instance, in \cite{Niculescu2004Stabilizing} in the case of a chain of integrators. Equation \eqref{eq:main-syst} corresponds to the particular case of \eqref{eq:general-LTI-delay} with $N = 2$, $\tau_0 = 0$, and $d = 1$, and its characteristic function $\Delta$ is given by
\begin{equation}
\label{eq:Delta}
\Delta(s) = s + a_0 + a_1 e^{-s \tau_1} + a_2 e^{-s \tau_2}.
\end{equation}
Note that $\Delta$ is a quasipolynomial of degree $3$, and hence any of its roots has multiplicity at most $3$. In other words, a root of \emph{maximal multiplicity} of $\Delta$ is a root of multiplicity $3$. Further insights on the stability regions in the delay-parameter space $(\tau_1,\tau_2)$ can be found in~\cite{hale-huang:93} (only scalar case) and \cite{gnc:05jmaa} (general two delays case). However, it is worth mentioning that the case of multiple characteristic roots is not explicitly addressed in any of the above cited references. Some ideas to deal with such a case can be found in \cite{freq-sweep-iter:19}, but the iterative frequency-sweeping methodology does not allow a deeper understanding of the existing links between maximal multiplicity and the location of the remaining characteristic roots. Our paper addresses this last issue.

\section{Definitions and preliminary results}
\label{sec:basic}

Since the goal of the paper is to study the MID property, which provides a link between multiplicity of a real root and its dominance, we start by providing a precise definition of dominance.

\begin{defn}
Let $\Delta: \mathbb C \to \mathbb C$ and $s_0 \in \mathbb R$ be such that $\Delta(s_0) = 0$. We say that $s_0$ is a \emph{dominant root} (resp.\ \emph{strictly dominant root}) of $\Delta$ if, for every $s \in \mathbb C$ such that $\Delta(s) = 0$, we have $\Real s \leq s_0$ (resp.\ $s = s_0$ or $\Real s < s_0$).
\end{defn}

The following proposition considers the effects of an affine transformation of the complex plane in the quasipolynomial $\Delta$ from \eqref{eq:Delta}.

\begin{prop}
\label{prop:normalization}
Let $a_0, a_1, a_2 \in \mathbb R$, $\tau_1, \tau_2 \in (0, +\infty)$ with $0 < \tau_1 < \tau_2$, $s_0 \in \mathbb R$, and consider the delay-differential equation \eqref{eq:main-syst} and its corresponding quasipolynomial $\Delta$ given by \eqref{eq:Delta}. Let $\widetilde\Delta: \mathbb C \to \mathbb C$ be defined for $z \in \mathbb C$ by $\widetilde\Delta(z) = \tau_2 \Delta(s_0 + \frac{z}{\tau_2})$. Then
\begin{equation}
\label{eq:TildeDelta}
\widetilde\Delta(z) = z + b_0 + b_1 e^{-\lambda z} + b_2 e^{-z},
\end{equation}
with
\begin{equation}
\label{eq:Relation-a-b}
\begin{aligned}
b_0 & = \tau_2(s_0 + a_0), & b_1 & = \tau_2 a_1 e^{-s_0 \tau_1}, \\
b_2 & = \tau_2 a_2 e^{-s_0 \tau_2}, & \lambda & = \frac{\tau_1}{\tau_2}.
\end{aligned}
\end{equation}
In particular, the function $s \mapsto \tau_2(s - s_0)$ is a bijection from roots of $\Delta$ to roots of $\widetilde\Delta$ preserving their multiplicities.
\end{prop}

Thanks to Proposition~\ref{prop:normalization}, a root $s_0$ of $\Delta$ is of maximal multiplicity, dominant, or strictly dominant if and only if the root $0$ of $\widetilde\Delta$ is of maximal multiplicity, dominant, or strictly dominant, respectively. Moreover, Proposition~\ref{prop:normalization} also shows that, with no loss of generality, one may replace the delays $(\tau_1, \tau_2)$ by $(\lambda, 1)$, i.e., one may assume that the largest delay is $1$.

\begin{prop}
\label{prop:max-multiplicity}
Let $a_0, a_1, a_2 \in \mathbb R$, $\tau_1, \tau_2 \in (0, +\infty)$ with $0 < \tau_1 < \tau_2$, $s_0 \in \mathbb R$, and consider the delay-differential equation \eqref{eq:main-syst} and its corresponding quasipolynomial $\Delta$ given by \eqref{eq:Delta}. Then $s_0$ is a root of $\Delta$ of maximal multiplicity $3$ if and only if
\begin{equation}
\label{eq:ConditionsA}
a_0 = -\frac{1}{\tau_1} - \frac{1}{\tau_2} - s_0, \, a_1 = \frac{\tau_2 e^{s_0 \tau_1}}{\tau_1 (\tau_2-\tau_1)}, \, a_2 = -\frac{\tau_1 e^{s_0 \tau_2}}{\tau_2(\tau_2-\tau_1)}.
\end{equation}
\end{prop}

\begin{proof}
Since $\Delta$ is a quasipolynomial of degree $3$, the maximal multiplicity of any of its roots is $3$. By Proposition~\ref{prop:normalization}, $s_0$ is a root of multiplicity $3$ of $\Delta$ if and only if $0$ is a root of multiplicity $3$ of the quasipolynomial $\widetilde\Delta$ defined by \eqref{eq:TildeDelta} and \eqref{eq:Relation-a-b}. Since $\widetilde\Delta$ is a quasipolynomial of degree $3$, $0$ is a root of multiplicity $3$ of $\widetilde\Delta$ if and only if $\widetilde\Delta(0) = \widetilde\Delta^\prime(0) = \widetilde\Delta^{\prime\prime}(0) = 0$. Thus, the conditions $b_0 + b_1 + b_2 = 0$, $1 - \lambda b_1 - b_2 = 0$, and $\lambda^2 b_1 + b_2 = 0$ should be simultaneously satisfied. The unique solution is
\begin{equation}
\label{eq:ConditionsB}
b_0 = -1 - \frac{1}{\lambda},\, b_1 = \frac{1}{\lambda (1 - \lambda)},\, b_2 = -\frac{\lambda}{1 - \lambda}.
\end{equation}
Using \eqref{eq:Relation-a-b}, one deduces that \eqref{eq:ConditionsB} is equivalent to \eqref{eq:ConditionsA}.
\end{proof}

Note that, under \eqref{eq:ConditionsB}, the quasipolynomial $\widetilde\Delta$ from \eqref{eq:TildeDelta} reads
\begin{equation}
\label{eq:TildeDeltaMID}
\widetilde\Delta(z) = z - 1 - \frac{1}{\lambda} + \frac{1}{\lambda (1 - \lambda)} e^{-\lambda z} - \frac{\lambda}{1 - \lambda} e^{-z}.
\end{equation}

\section{A rationally dependent case}
\label{sec:1/2}

One of the classical arguments to prove the MID property for single-delay systems, used, for instance, in \cite{Boussaada_ESAIM}, is to make use of Cauchy's argument principle to show that no roots lie to the right of the root with maximal multiplicity. In this section, we apply this technique to \eqref{eq:main-syst} in the case $\frac{\tau_1}{\tau_2} = \frac{1}{2}$, in which several simplifications can be carried out due to the simple rational dependence of the delays. We rely on the results of \cite{Stepan1979Stability} and \cite{hassard1997counting}, in which the authors provide formulas for the number of roots of a quasipolynomial on the right half-plane after an application of Cauchy's argument principle.

\begin{thm}
Let $a_0, a_1, a_2 \in \mathbb R$, $\tau_1, \tau_2 \in (0, +\infty)$ with $0 < \tau_1 < \tau_2$, $s_0 \in \mathbb R$, consider the delay-differential equation \eqref{eq:main-syst} and its corresponding quasipolynomial $\Delta$ given by \eqref{eq:Delta}, and assume that $\frac{\tau_1}{\tau_2} = \frac{1}{2}$. If $s_0$ is a root of maximal multiplicity $3$ of $\Delta$, then $s_0$ is a strictly dominant root of $\Delta$. In particular, if $s_0 < 0$, then \eqref{eq:main-syst} is exponentially stable.
\end{thm}

\begin{proof}
Let $\widetilde\Delta$ be obtained from $\Delta$ as in Proposition~\ref{prop:normalization}. If $s_0$ is a root of maximal multiplicity of $\Delta$ and $\lambda = \frac{\tau_1}{\tau_2} = \frac{1}{2}$, then, by \eqref{eq:TildeDeltaMID},
\[
\widetilde\Delta(z) = z - 3 + 4 e^{-\frac{z}{2}} - e^{-z}.
\]
Moreover, by Proposition~\ref{prop:normalization}, $s_0$ is a strictly dominant root of $\Delta$ if and only if $0$ is a strictly dominant root of $\widetilde\Delta$.

Let us then show that $0$ is a strictly dominant root of $\widetilde\Delta$. For that purpose, we will apply the main result of \cite{hassard1997counting}, which provides an explicit expression for the number of zeros of a quasipolynomial in the right half-plane. As a preliminary step, we first show that the unique root of $\widetilde\Delta$ on the imaginary axis is $0$.

\begin{claim}
\label{claim:noRootOnImaginaryAxis}
If $y \in \mathbb R\setminus \{0\}$, then $\widetilde\Delta(i y) \neq 0$.
\end{claim}

\begin{proof*}{\textnormal{\emph{Proof of Claim~\ref{claim:noRootOnImaginaryAxis}.}}}
Let $y \in \mathbb R$ and assume that $\widetilde\Delta(i y) = 0$, i.e., $iy - 3 + 4 e^{-i\frac{y}{2}} - e^{-iy} = 0$. In particular, $\abs{3 - i y} = \abs{4 e^{-i\frac{y}{2}} - e^{-iy}}$ and thus $9 + y^2 \leq 25$, implying that $\abs{y} \leq 4$. Letting $\zeta = \frac{y}{2}$, we then have that $\abs{\zeta} \leq 2$ and $2 i \zeta - 3 + 4 e^{- i \zeta} - e^{- 2 i \zeta} = 0$. Taking real and imaginary parts and using standard trigonometric identities, we deduce that
\begin{equation*}
(2 - \cos\zeta) \cos\zeta = 1, \qquad
(2 - \cos\zeta) \sin\zeta = \zeta,
\end{equation*}
which implies that $\tan\zeta = \zeta$. Since the unique solution of the latter equation in the interval $[-2, 2]$ is $\zeta = 0$, we finally obtain that $y = 0$, as required.
\end{proof*}

In order to apply the main result of \cite{hassard1997counting}, let us introduce some notation. Let $Z$ be the number of roots of $\widetilde\Delta$ with positive real part counted by multiplicity, $R: \mathbb R \to \mathbb R$ and $S: \mathbb R \to \mathbb R$ be defined for $y \in \mathbb R$ by $R(y) = -\Real[i \widetilde\Delta(i y)]$ and $S(y) = -\Imag[i \widetilde\Delta(i y)]$, denote by $r \in \mathbb N$ the number of zeros of $R$ in $(0, +\infty)$, and let $\rho_1, \dotsc, \rho_r$ be the positive zeros of $R$, repeated according to their multiplicity and ordered so that $0 < \rho_1 \leq \dotsb \leq \rho_r$. Using Claim~\ref{claim:noRootOnImaginaryAxis}, the main result of \cite{hassard1997counting} implies that
\begin{equation}
\label{eq:Hassard}
Z = -1 + \frac{(-1)^r}{2} \sign S^{(3)}(0) + \sum_{j=1}^r (-1)^{j-1} \sign S(\rho_j),
\end{equation}
where $\sign: \mathbb R \to \mathbb R$ is the sign function, defined by $\sign(0) = 0$ and $\sign(x) = \frac{x}{\abs{x}}$ for $x \neq 0$.

We have $R(y) = y - 4 \sin(y/2) + \sin(y)$. For $y \geq 2 \pi$, we have $R(y) > 0$, and then all positive zeros of $R$ belong to the interval $(0, 2\pi)$. We have $R^\prime(y) = 2 \cos(y/2) \left[\cos(y/2) - 1\right]$, and thus $R$ is decreasing in $[0, \pi]$ and increasing in $[\pi, 2 \pi]$. Since $R(0) = 0$, $R(\pi) = \pi - 4 < 0$, and $R(2\pi) = 2\pi > 0$, we deduce that $R$ admits a unique positive zero $\rho_1$, which belongs to $(\pi, 2\pi)$. Then $r = 1$ and $\rho_1 \in (\pi, 2 \pi)$, and thus \eqref{eq:Hassard} becomes
\[
Z = -1 - \frac{1}{2} \sign S^{(3)}(0) + \sign S(\rho_1).
\]
Using the fact that $Z$ is a nonnegative integer and that $\sign(x) \in \{-1, 0, 1\}$ for every $x \in \mathbb R$, one deduces at once that $\sign S^{(3)}(0) = 0$, $\sign S(\rho_1) = 1$, and $Z = 0$. Hence $\widetilde\Delta$ admits no roots at the open right half-plane and, combining with Claim~\ref{claim:noRootOnImaginaryAxis}, we deduce that $0$ is a strictly dominant root of $\widetilde\Delta$, as required.
\end{proof}

The above proof can be adapted to some other ``simple'' rationally dependent cases, such as $\frac{\tau_1}{\tau_2} = \frac{1}{3}$, $\frac{2}{3}$, $\frac{1}{4}$, $\frac{3}{4}$, etc., but a general formulation for every rational $\frac{\tau_1}{\tau_2}$ seems out of reach, the major difficulties being the generalization of Claim~\ref{claim:noRootOnImaginaryAxis} and the analysis of the behavior of $R$ and $S$, which becomes much more delicate. For this reason, we consider, in the next section, another strategy for proving the MID property for \eqref{eq:main-syst} in a more general setting.

\section{Bounded delay ratio}
\label{sec:max-ratio}

In this section, we consider another approach to prove the MID property for \eqref{eq:main-syst}, inspired by the Walton--Marshall method from \cite{Walton1987Direct} and similar to the strategy used in \cite{Mazanti2020Spectral} to study the MID property for complex conjugate dominant roots. This approach consists on regarding the quasipolynomial $\widetilde\Delta$ from \eqref{eq:TildeDeltaMID} as a quasipolynomial depending on a parameter $\lambda \in (0, 1)$. Standard arguments allow one to prove that any root of $\widetilde\Delta$ is a continuous function of $\lambda$ (and actually analytic in neighborhoods of values of $\lambda$ at which this root is simple), which can be extended at $\lambda = 0$. The strategy is then to prove that $0$ is a dominant root in the limit case $\lambda = 0$ and that this dominance is preserved as $\lambda$ increases, by showing that no root may cross the imaginary axis and that roots coming from infinity cannot appear in the right half-plane.

For that purpose, we find it useful to normalize the quasipolynomial $\widetilde\Delta$ from \eqref{eq:TildeDeltaMID} by multiplying it by $\frac{1}{\lambda}$ and to consider it as a function of both variables $z$ and $\lambda$. Hence, we consider in the sequel the function $\widehat\Delta: \mathbb C \times (0, 1) \to \mathbb C$ defined by
\begin{equation}
\label{eq:HatDelta}
\widehat\Delta(z, \lambda) = \frac{z}{\lambda} - \frac{1}{\lambda} - \frac{1}{\lambda^2} + \frac{1}{\lambda^2 (1 - \lambda)} e^{-\lambda z} - \frac{1}{1 - \lambda} e^{-z}.
\end{equation}
Several useful properties of $\widehat\Delta$ are presented in Appendix~\ref{app:HatDelta}, in particular the fact that $\widehat\Delta$ can be extended in a unique way to a holomorphic function defined in $\mathbb C^2$.

\begin{thm}
\label{thm:max-ratio}
Let $a_0, a_1, a_2 \in \mathbb R$, $\tau_1, \tau_2 \in (0, +\infty)$ with $0 < \tau_1 < \tau_2$, $s_0 \in \mathbb R$, consider the delay-differential equation \eqref{eq:main-syst} and its corresponding quasipolynomial $\Delta$ given by \eqref{eq:Delta}, and assume that $\frac{\tau_1}{\tau_2} \in (0, \frac{2}{3}]$. If $s_0$ is a root of maximal multiplicity $3$ of $\Delta$, then $s_0$ is a strictly dominant root of $\Delta$. In particular, if $s_0 < 0$, then \eqref{eq:main-syst} is exponentially stable.
\end{thm}

\begin{proof}
Thanks to Proposition~\ref{prop:normalization}, it suffices to show that, for every $\lambda \in [0, 1]$, the root of maximal multiplicity $0$ of the quasipolynomial $\widehat\Delta(\cdot, \lambda)$ is strictly dominant. By Lemma~\ref{lem:MIDHatDelta0} in Appendix~\ref{app:HatDelta}, this is the case for $\lambda = 0$. Assume, to obtain a contradiction, that there exist $\lambda_\ast \in (0, \frac{2}{3}]$ and $z_\ast \in \mathbb C \setminus \{0\}$ such that $\Real z_\ast \geq 0$ and $\widehat\Delta(z_\ast, \lambda_\ast) = 0$. Thanks to Lemma~\ref{lem:unique-on-imaginary} in Appendix~\ref{app:HatDelta}, one has necessarily $\Real z_\ast > 0$.

Since $\widehat\Delta$ is holomorphic, there exist an interval $\Lambda \subset (0, \frac{2}{3}]$ containing $\lambda_\ast$ and a continuous function $\zeta: \Lambda \to \mathbb C$ such that $\zeta(\lambda_\ast) = z_\ast$ and $\widehat\Delta(\zeta(\lambda), \lambda) = 0$ for every $\lambda \in \Lambda$. Note that, since $0$ is a root of constant multiplicity $3$ of $\widehat\Delta(\cdot, \lambda)$ for every $\lambda \in (0, \frac{2}{3}]$, Hurwitz Theorem ensures that $\zeta(\lambda) \neq 0$ for every $\lambda \in \Lambda$. By Lemma~\ref{lem:unique-on-imaginary} in Appendix~\ref{app:HatDelta}, we then deduce that $\Real \zeta(\lambda) \neq 0$ for every $\lambda \in \Lambda$, and thus, since $\Real \zeta(\lambda_\ast) > 0$ and $\zeta$ is continuous, we have that $\Real \zeta(\lambda) > 0$ for every $\lambda \in \Lambda$. Using this fact and \eqref{eq:HatDelta}, we get the bound
\begin{equation}
\label{eq:bound-zeta-lambda}
\abs{\zeta(\lambda)} \leq 1 + \frac{1}{\lambda} + \frac{1}{\lambda(1 - \lambda)} + \frac{\lambda}{1 - \lambda} = \frac{2}{\lambda(1 - \lambda)}
\end{equation}
for every $\lambda \in \Lambda$, which shows that $\zeta$ can be extended to a continuous function (still denoted by $\zeta$) defined on $(0, \frac{2}{3}]$, which still satisfies $\Real \zeta(\lambda) > 0$ for every $\lambda \in (0, \frac{2}{3}]$.

Let us now show that $\zeta(\lambda)$ remains bounded as $\lambda \to 0^+$. Using \eqref{eq:HatDelta}, we have
\begin{equation*}
\lambda \zeta(\lambda) = \frac{1 - e^{-\lambda \zeta(\lambda)}}{1-\lambda} - \frac{\lambda^2 (1-e^{-\zeta(\lambda)})}{1-\lambda}.
\end{equation*}
In particular, $\lambda \mapsto \lambda \zeta(\lambda)$ is bounded on $(0, \frac{2}{3}]$. Let $\xi(\lambda) = \lambda \zeta(\lambda)$ and $g: \mathbb C \to \mathbb C$ be the entire function defined for $z \in \mathbb C \setminus \{0\}$ by $g(z) = \frac{1-e^{-z}}{z} - 1$. Then
\begin{equation*}
\xi(\lambda) = \frac{\xi(\lambda)\left[1 + g(\xi(\lambda))\right]}{1-\lambda} + O(\lambda^2)
\end{equation*}
as $\lambda \to 0^+$. Then
\begin{equation}
\label{eq:lem_en_0}
-\lambda \xi(\lambda) = \xi(\lambda) g(\xi(\lambda)) + O(\lambda^2).
\end{equation}
Since $\xi$ is bounded, the left-hand side of \eqref{eq:lem_en_0} converges to $0$ as $\lambda \to 0^+$, and thus $\xi(\lambda) g(\xi(\lambda))$ also converges to zero as $\lambda \to 0^+$. Since the only zero of the equation $z g(z) = 0$ in the half-space $\{z \in \mathbb C \mid \Real z \ge 0\}$ is $z = 0$, we deduce that $\xi(\lambda) \to 0$ as $\lambda \to 0^+$. We have $g(z) = -z(1+o(1))/2$ as $z \to 0$ and, inserting this into \eqref{eq:lem_en_0}, we obtain that
\begin{equation*}
-\lambda \xi(\lambda) = -\frac{\xi(\lambda)^2}{2}(1+o(1)) + O(\lambda^2).
\end{equation*}
Then $\xi(\lambda) = O(\lambda)$ and $\zeta(\lambda) = O(1)$ as $\lambda \to 0^+$, concluding the proof that $\zeta(\lambda)$ remains bounded as $\lambda \to 0^+$.

Hence, there exist $z_0 \in \mathbb C$ and a sequence $(\lambda_n)_{n \in \mathbb N}$ in $(0, \frac{2}{3}]$ such that $\lambda_n \to 0^+$ and $\zeta(\lambda_n) \to z_0$ as $n \to +\infty$. Then, by Hurwitz Theorem, $z_0$ is a root of $\widehat\Delta(\cdot, 0)$ and, since $\Real \zeta(\lambda) > 0$ for every $\lambda \in (0, \frac{2}{3}]$, we have that $\Real z_0 \geq 0$. From Lemma~\ref{lem:MIDHatDelta0} in Appendix~\ref{app:HatDelta}, we conclude that $z_0 = 0$. However, since $0$ is a root of multiplicity $3$ of $\widehat\Delta(\cdot, \lambda)$ for every $\lambda \in [0, \frac{2}{3}]$ and $\zeta(\lambda) \neq 0$ for $\lambda \in (0, \frac{2}{3}]$, Hurwitz Theorem also implies that $0$ should be a root of multiplicity at least $4$ of $\widehat\Delta(\cdot, 0)$, which contradicts Lemma~\ref{lem:MIDHatDelta0}. This contradiction establishes the result.
\end{proof}

The major difficulty in extending the above technique to any $\frac{\tau_1}{\tau_2} \in (0, 1)$ comes from extending Lemma~\ref{lem:unique-on-imaginary} in Appendix~\ref{app:HatDelta}. Indeed, the method used in that lemma relies on the analysis of the function $F$ defined in \eqref{eq:F}. The result of Lemma~\ref{lem:unique-on-imaginary} actually holds for $\lambda \in (0, \frac{1}{F_{\max}}) \cup (F_{\max}, +\infty)$, where $F_{\max} = \max_{\omega > 0} F(\omega)$, and its proof shows that $F_{\max} \leq \frac{4\pi^2}{(2\pi - 1)^2}$, yielding that the statement of Theorem~\ref{thm:max-ratio} still holds with $(0, \frac{2}{3}]$ replaced by $(0, (1 - \frac{1}{2\pi})^2)$. (Note that $(1 - \frac{1}{2\pi})^2 \approx 0.707$.) Numerically, we have $F_{\max} \approx 1.264$, and thus, using the same strategy, one can expect to obtain the statement of Theorem~\ref{thm:max-ratio} for $\lambda \in (0, \frac{1}{F_{\max}})$, where $\frac{1}{F_{\max}} \approx 0.791$.

\section{Numerical illustration}
\label{sec:numerical}

In view of Theorem~\ref{thm:max-ratio}, a natural question is whether the MID property holds for \eqref{eq:main-syst} for every $\tau_1, \tau_2 \in (0, +\infty)$ with $0 < \tau_1 < \tau_2$, i.e., whether $0$ is a strictly dominant root for the quasipolynomial $\widehat\Delta(\cdot, \lambda)$ from \eqref{eq:HatDelta} for every $\lambda \in (0, 1)$.

\begin{figure}[ht]
\centering
\resizebox{0.95\columnwidth}{!}{\input{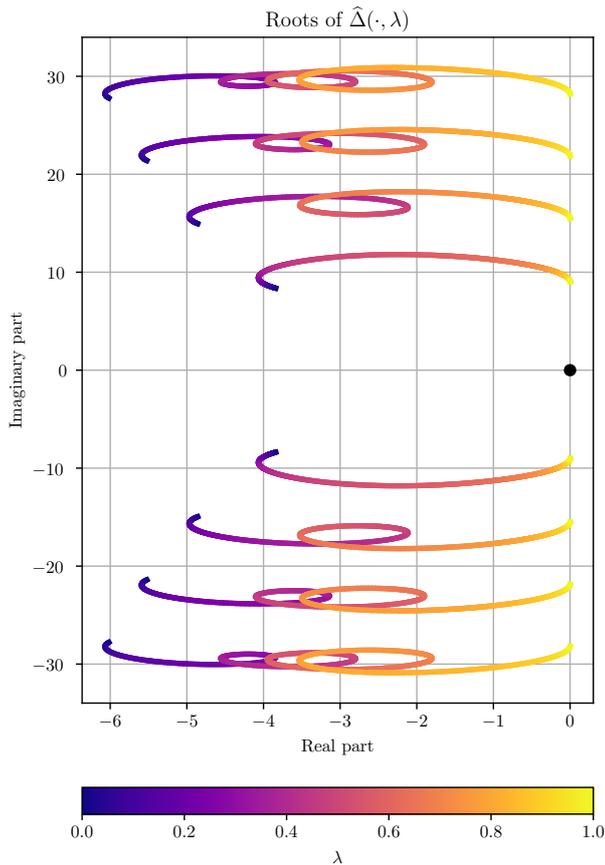}}
\caption{Roots of $\widehat\Delta(\cdot, \lambda)$ for $\lambda \in [0, 1]$.}
\label{FigSpectrum}
\end{figure}

We represent, in Figure~\ref{FigSpectrum}, the numerical computation of the roots of $\widehat\Delta(\cdot, \lambda)$ as $\lambda$ goes from $0$ to $1$ (color scale from blue to yellow), with the black dot at the origin representing the root $0$ of multiplicity $3$ common to $\widehat\Delta(\cdot, \lambda)$ for every $\lambda \in [0, 1]$. Only roots in the rectangle $\{z \in \mathbb C \mid -10 \leq \Real z \leq 5,\, -33 \leq \Imag z \leq 33\}$ were computed, and all numerical computations were done using Python's \texttt{cxroots} package from \cite{cxroots}, which is based on methods from \cite{Kravanja2000Computing}.

We observe in Figure~\ref{FigSpectrum} that, as shown in Theorem~\ref{thm:max-ratio}, all roots of $\widehat\Delta(\cdot, \lambda)$ represented in the figure are strictly to the left of the root of maximal multiplicity $0$ for every $\lambda \in (0, \frac{2}{3}]$, and actually for every $\lambda \in [0, 1)$. For $\lambda = 1$, as previously shown in \cite{MBNC}, all roots lie on the imaginary axis.

\section{General statement}
\label{sec:general}

It turns out that one can refine the arguments from the proof of Theorem~\ref{thm:max-ratio} to show that its statement holds for every $\tau_1, \tau_2 \in (0, +\infty)$ with $0 < \tau_1 < \tau_2$, without any additional restriction on $\frac{\tau_1}{\tau_2}$. More precisely, we have the following result.

\begin{thm}
\label{thm:general}
Let $a_0, a_1, a_2 \in \mathbb R$, $\tau_1, \tau_2 \in (0, +\infty)$ with $0 < \tau_1 < \tau_2$, $s_0 \in \mathbb R$, and consider the delay-differential equation \eqref{eq:main-syst} and its corresponding quasipolynomial $\Delta$ given by \eqref{eq:Delta}. If $s_0$ is a root of maximal multiplicity $3$ of $\Delta$, then $s_0$ is a strictly dominant root of $\Delta$. In particular, if $s_0 < 0$, then \eqref{eq:main-syst} is exponentially stable.
\end{thm}

Let us briefly sketch the strategy proof of Theorem~\ref{thm:general}, which will be detailed in an upcoming paper. We follow the same strategy of the proof of Theorem~\ref{thm:max-ratio}, noticing first that, by Lemma~\ref{lem:MIDHatDelta0}, the root of maximal multiplicity $0$ is strictly dominant for $\widehat\Delta(\cdot, 0)$ and then showing that this dominance is preserved as $\lambda$ increases. To do that, as in the proof of Theorem~\ref{thm:max-ratio}, we assume, in order to obtain a contradiction, that there exist $\lambda_0 \in (0, 1)$ and $z_0 \in \mathbb C \setminus \{0\}$ with $\Real z_0 \geq 0$ such that $\widehat\Delta(z_0, \lambda_0) = 0$. Such a root $z_0$ of $\widehat\Delta(\cdot, \lambda_0)$ is necessarily simple and hence there exists an analytic function $\zeta$ defined on a neighborhood $\Lambda$ of $\lambda_0$ such that $\zeta(\lambda_0) = z_0$ and $\widehat\Delta(\zeta(\lambda), \lambda) = 0$ for every $\lambda \in \Lambda$. Using \eqref{eq:bound-zeta-lambda}, $\zeta$ can be extended to a continuous function defined on $(0, 1)$ still satisfying the same properties and, as in the proof of Theorem~\ref{thm:max-ratio}, this function is bounded in a neighborhood of $0$. Using the fact that $0$ is a strictly dominant roof of $\widehat\Delta(\cdot, 0)$ and arguing as in the proof of Theorem~\ref{thm:max-ratio}, we deduce that there exists $\lambda_1 \in (0, \lambda_0]$ such that $\Real \zeta(\lambda_1) = 0$, i.e., the root $\zeta(\lambda_1)$ of $\widehat\Delta(\cdot, \lambda_1)$ lies on the imaginary axis.

The next step is then to show that, if $\zeta(\lambda)$ lies on the imaginary axis for some $\lambda \in (0, 1)$, then necessarily $\zeta(\lambda)$ is a simple root of $\widehat\Delta(\cdot, \lambda)$ and $\Real \zeta^\prime(\lambda) > 0$. This will show, in particular, that $\Real \zeta(\lambda) > 0$ for every $\lambda \in (\lambda_1, 1)$. We can then show that $\zeta$ remains bounded for $\lambda$ close to $1$, proving that $\zeta(\lambda)$ must converge, as $\lambda \to 1^-$, to one of the roots of $\widehat\Delta(\cdot, 1)$ on the imaginary axis apart from the origin. However, on those roots, we compute $\Real \zeta^\prime(1) = 0$ and $\Real \zeta^{\prime\prime}(1) < 0$, yielding that $\Real \zeta(\lambda) < 0$ for $\lambda$ close to $1$, and hence contradicting the fact that $\Real \zeta(\lambda) > 0$ for every $\lambda \in (\lambda_1, 1)$. This contradiction finally yields the conclusion.

\bibliography{ifacconf}

\appendix

\section{Properties of $\widehat\Delta$}
\label{app:HatDelta}

In this section, we provide several properties of the function $\widehat\Delta$ defined in \eqref{eq:HatDelta} that are useful for the proof of Theorem~\ref{thm:max-ratio} in Section~\ref{sec:max-ratio}.

\begin{lem}
Let $\widehat\Delta: \mathbb C \times (0, 1) \to \mathbb C$ be defined by \eqref{eq:HatDelta}. Then $\widehat\Delta$ can be extended in a unique manner to a holomorphic function (still denoted by $\widehat\Delta$) defined on $\mathbb C^2$. Moreover, this extension satisfies
\begin{subequations}
\label{eq:HatDelta-Limit}
\begin{align}
\widehat\Delta(z, 0) & = \frac{z^2}{2} - z + 1 - e^{-z}, \label{eq:HatDelta-0} \\
\widehat\Delta(z, 1) & = z - 2 + (z + 2) e^{-z}. \label{eq:HatDelta-1}
\end{align}
\end{subequations}
\end{lem}

\begin{proof}
Clearly, $\widehat\Delta$ can be extended to a holomorphic function defined on $\mathbb C \times (\mathbb C \setminus \{0, 1\})$ and, to conclude the proof, it suffices to show that $\lambda = 0$ and $\lambda = 1$ are actually removable singularities of $\widehat\Delta$. Hence, the proof is completed if one shows that $\lim_{\lambda \to 0} \widehat\Delta(z, \lambda)$ and $\lim_{\lambda \to 1} \widehat\Delta(z, \lambda)$ exist and are given by the respective expressions in \eqref{eq:HatDelta-Limit}.

As $\lambda \to 0$, we have 
\begin{multline*}
\lambda^2 \widehat\Delta(z, \lambda) = \lambda z - \lambda - 1 \\
{} + \left(1 - \lambda z + \frac{\lambda^2 z^2}{2} + O(\lambda^3)\right) (1 + \lambda + \lambda^2 + O(\lambda^3)) \\
{} - \lambda^2 e^{-z} (1 + O(\lambda))\\
= \lambda^2 \left[\frac{z^2}{2} - z + 1 - e^{-z} + O(\lambda)\right],
\end{multline*}
uniformly with respect to $z$ on compact subsets of $\mathbb C$. Hence, letting $\lambda \to 0$, we deduce that the singularity of $\widehat\Delta$ at $\lambda = 0$ is removable and that \eqref{eq:HatDelta-0} holds.

Letting now $\lambda \to 1$ and setting $\beta = 1 - \lambda$ for convenience, we write the last two terms of \eqref{eq:HatDelta} as
\begin{multline*}
\frac{\frac{1}{\lambda^2} e^{-\lambda z} - e^{-z}}{1 - \lambda} = \frac{e^{-z}}{\beta} \left[\frac{1}{(1 - \beta)^2} e^{\beta z} - 1\right]\\ = \frac{e^{-z}}{\beta} \left[(1 + 2\beta + O(\beta^2)) (1 + \beta z + O(\beta^2)) - 1\right] \\ = e^{-z} \left[z + 2 + O(1 - \lambda)\right],
\end{multline*}
uniformly with respect to $z$ on compact subsets of $\mathbb C$. Inserting into \eqref{eq:HatDelta} and letting $\lambda \to 1$ shows that the singularity of $\widehat\Delta$ at $\lambda = 1$ is removable and thus \eqref{eq:HatDelta-1} holds.
\end{proof}

\begin{rem}
The quasipolynomial $\widehat\Delta(\cdot, 0)$ is the characteristic function of the second-order retarded delay-differential equation with a single delay
\[
\frac{1}{2} y^{\prime\prime}(t) - y^\prime(t) + y(t) - y(t - 1) = 0,
\]
whereas $\widehat\Delta(\cdot, 1)$ is the characteristic function of the first-order neutral delay-differential equation with a single delay
\[
y^\prime(t) - 2 y(t) + y^\prime(t - 1) + 2 y(t - 1) = 0.
\]
The roots of $\widehat\Delta(\cdot, 1)$ have been completely characterized on \cite{MBNC} and, in particular, they all lie on the imaginary axis.
\end{rem}

Let us now state the following symmetry property of $\widehat\Delta$, whose proof is immediate.

\begin{lem}
\label{lem:symmetry}
Let $\widehat\Delta$ be defined on $\mathbb C^2$ by \eqref{eq:HatDelta}. Then, for every $z \in \mathbb C$ and $\lambda \in \mathbb C \setminus \{0\}$, we have
\[
\lambda^3 \widehat\Delta(z, \lambda) = \widehat\Delta\left(\lambda z, \frac{1}{\lambda}\right).
\]
\end{lem}

We next show, using a factorization technique, that the multiplicity-induced-dominancy property holds for $\widehat\Delta(\cdot, 0)$.

\begin{lem}
\label{lem:MIDHatDelta0}
Consider the quasipolynomial $\widehat\Delta(\cdot, 0)$ defined in \eqref{eq:HatDelta-0}. Then $0$ is a root of multiplicity $3$ of $\widehat\Delta(\cdot, 0)$ and it is strictly dominant.
\end{lem}

\begin{proof}
It is immediate to verify from \eqref{eq:HatDelta-0} that $0$ is a root of multiplicity $3$ of $\widehat\Delta(\cdot, 0)$. Moreover, an explicit computation shows that $\widehat\Delta(\cdot, 0)$ can be factorized, for every $z \in \mathbb C$, as
\begin{equation}
\label{eq:factorization-0}
\widehat\Delta(z, 0) = z^2\left(\frac{1}{2} - \int_0^1 \int_0^t e^{-z \xi} \diff \xi \diff t\right).
\end{equation}
If $z \in \mathbb C$ is such that $\Real z > 0$, then $\abs*{\int_0^1 \int_0^t e^{-z \xi} \diff \xi \diff t} \leq \int_0^1 \int_0^t e^{-(\Real z) \xi} \diff \xi \diff t < \frac{1}{2}$ since $e^{-(\Real z) \xi} < 1$ for $\xi > 0$. The factorization \eqref{eq:factorization-0} then implies that such a $z$ cannot be a root of $\widehat\Delta(\cdot, 0)$, and thus $0$ is a dominant root of $\widehat\Delta(\cdot, 0)$. To show that it is strictly dominant, we are left to show that the only root of $\widehat\Delta(\cdot, 0)$ lying on the imaginary axis is $0$. Indeed, let $y \in \mathbb R$ be such that $\widehat\Delta(iy, 0) = 0$. Then $-\frac{y^2}{2} - iy + 1 = e^{-i y}$ and, taking the imaginary part, we deduce that $y = \sin y$, which implies that $y = 0$, as required.
\end{proof}

The last useful property of $\widehat\Delta$ shown in this appendix is that, for a range of real values of $\lambda$ containing those considered in Theorem~\ref{thm:max-ratio}, $\widehat\Delta(\cdot, \lambda)$ admits no roots on the imaginary axis apart from $0$.

\begin{lem}
\label{lem:unique-on-imaginary}
Consider the quasipolynomial $\widehat\Delta$ from \eqref{eq:HatDelta} and assume that $\lambda \in (0, \frac{2}{3}] \cup [\frac{3}{2}, +\infty)$. Then the unique root of $\widehat\Delta(\cdot, \lambda)$ on the imaginary axis is $0$.
\end{lem}

\begin{proof}
Let $F: (0, +\infty) \to \mathbb R$ be defined by
\begin{equation}
\label{eq:F}
F(\omega) = \frac{\omega^2 + 2(\cos\omega-1)}{(\omega-\sin\omega)^2+(1-\cos\omega)^2}.
\end{equation}

\begin{claim}
\label{claim:nec-condition-root}
If $\omega, \mu \in \mathbb R$ are such that $\omega \neq 0$, $\mu > 0$, $\mu \neq 1$, and $\widehat\Delta(i \omega, \mu) = 0$, then $\mu = F(\omega)$.
\end{claim}

\begin{proof*}{\textnormal{\emph{Proof of Claim~\ref{claim:nec-condition-root}.}}}
Indeed, from \eqref{eq:HatDelta}, $\widehat\Delta(i \omega, \mu) = 0$ reads
\[
i \omega - 1 - \frac{1}{\mu} - \frac{\mu}{1 - \mu} e^{-i \omega} = - \frac{e^{- i \mu \omega}}{\mu (1 - \mu)}.
\]
Thus
\[
\abs*{i \omega - 1 - \frac{1}{\mu} - \frac{\mu}{1 - \mu} e^{-i \omega}}^2 = \frac{1}{\mu^2 (1 - \mu)^2}
\]
and, developing this expression, we get that
\begin{equation*}
2(1 + \mu)(1-\cos\omega) - (1-\mu) \omega^2 - 2 \mu \omega \sin\omega = 0,
\end{equation*}
which implies that $\mu = F(\omega)$, as required.
\end{proof*}

\begin{claim}
\label{claim:bound-F}
The function $F$ is odd and $F(\omega) < \frac{3}{2}$ for every $\omega \in \mathbb R \setminus \{0\}$.
\end{claim}

\begin{proof*}{\textnormal{\emph{Proof of Claim~\ref{claim:bound-F}.}}}
The fact that $F$ is odd is immediate. Let us first study the behavior of $F$ on $(0, 2\pi]$. We claim that, in this interval, we have $F(\omega) \leq 1$. Indeed, notice that $F(\omega) \leq 1$ is equivalent to $2 (1 - \cos\omega) \geq \omega \sin \omega$. Defining $G: \mathbb R \to \mathbb R$ by $G(\omega) = 2(1 - \cos\omega) - \omega \sin\omega$, we have that $G(0) = G(2 \pi) = 0$ and $G^\prime(\omega) = \sin\omega - \omega \cos\omega$. Then $G$ is increasing on $[0, \omega_\ast]$ and decreasing on $[\omega_\ast, 2\pi]$, where $\omega_\ast$ is the smallest positive solution of $\tan \omega_\ast = \omega_\ast$, and thus $G(\omega) \geq 0$ for every $\omega \in [0, 2 \pi]$, implying that $F(\omega) \leq 1$ for every $\omega \in (0, 2\pi]$. Finally, if $\omega > 2 \pi$, we estimate $F$ as $F(\omega) \leq \frac{\omega^2}{(\omega-1)^2} \leq \frac{4 \pi^2}{(2\pi - 1)^2} < \frac{3}{2}$.
\end{proof*}

Now, let $y \in \mathbb R$ and $\lambda \in (0, \frac{2}{3}] \cup [\frac{3}{2}, +\infty)$ be such that $\widehat\Delta(i y, \lambda) = 0$ and assume, to obtain a contradiction, that $y \neq 0$. By Lemma~\ref{lem:symmetry}, we also have $\widehat\Delta\left(i \lambda y, \frac{1}{\lambda}\right) = 0$, and thus, by Claim~\ref{claim:nec-condition-root}, we deduce that we have both $\lambda = F(y)$ and $\frac{1}{\lambda} = F(\lambda y)$. Hence, by Claim~\ref{claim:bound-F}, we have that $\lambda < \frac{3}{2}$ and $\frac{1}{\lambda} < \frac{3}{2}$, which contradicts the fact that $\lambda \in (0, \frac{2}{3}] \cup [\frac{3}{2}, +\infty)$. This contradiction establishes the result.
\end{proof}

\end{document}